# ESTIMATING INVARIANT LAWS OF LINEAR PROCESSES BY U-STATISTICS

By Anton Schick[1] and Wolfgang Wefelmeyer

*Binghamton University and Universität zu Köln*

Suppose we observe an invertible linear process with independent mean-zero innovations and with coefficients depending on a finite-dimensional parameter, and we want to estimate the expectation of some function under the stationary distribution of the process. The usual estimator would be the empirical estimator. It can be improved using the fact that the innovations are centered. We construct an even better estimator using the representation of the observations as infinite-order moving averages of the innovations. Then the expectation of the function under the stationary distribution can be written as the expectation under the distribution of an infinite series in terms of the innovations, and it can be estimated by a $U$-statistic of increasing order (also called an "infinite-order $U$-statistic") in terms of the estimated innovations. The estimator can be further improved using the fact that the innovations are centered. This improved estimator is optimal if the coefficients of the linear process are estimated optimally. The variance reduction of our estimator over the empirical estimator can be considerable.

**1. Introduction.** There is a large literature on estimation in ergodic time series driven by independent innovations. In the last fifteen years, optimality questions have also been addressed. Efficient estimators for the parameters of ARMA-type processes are constructed by Kreiss (1987a, b), Jeganathan (1995), Drost, Klaassen and Werker (1997), Koul and Schick (1997) and Schick and Wefelmeyer (2002a). For invertible linear time series, the innovations can be estimated, and linear functionals of the innovation distribution can then be estimated by corresponding empirical estimators based on the estimated innovations; see Boldin (1982) and Kreiss (1991). Simple

Received August 2001; revised April 2003.
[1]Supported in part by NSF Grant DMS 00-72174.

and efficient improvements of these estimators are possible if the innovations are centered; see Wefelmeyer (1994) and Schick and Wefelmeyer (2002b).

Here we are interested in estimating functionals of the stationary law. Such functionals can be estimated in a straightforward way from observations of the time series. Linear functionals of the stationary law can be estimated by corresponding empirical estimators. The stationary density can be estimated by a kernel estimator; see, for example, Yakowitz (1989), Tran (1992) and Honda (2000).

These estimators are "nonparametric" in that they do not exploit the information that the time series is driven by independent innovations. In this paper we show how to use this information in order to construct efficient estimators for *linear functionals* of the stationary law of causal and invertible linear processes with coefficients depending on a finite-dimensional parameter. We restrict attention to estimation of expectations of smooth functions. Examples are moments, absolute moments, the characteristic function and other transformations of the stationary law. One of the applications would be testing for Gaussianity. Under stronger conditions on the time series, one could prove corresponding results for expectations of step functions, for example, the distribution function. An application would be estimating the value at risk in financial mathematics.

In the simplest such time series, a moving average process of order 1, Saavedra and Cao (1999, 2000) show that the specific structure of the model allows the stationary *density* to be estimated at the parametric rate $n^{-1/2}$. Schick and Wefelmeyer (2004) prove that the estimator of Saavedra and Cao is efficient. Analogous parametric rates can also be obtained for estimators of *conditional expectations*; see Müller, Schick and Wefelmeyer (2003) for a result in nonlinear autoregressive processes. Such estimators could be combined with the estimators in the present paper in order to efficiently estimate functionals of *joint laws* of linear processes, for example, autocovariance functions.

A cautionary remark: unlike the usual empirical estimators for functionals of the stationary law, our efficient estimators use the full structure of the model, in particular, the independence of the innovations. Like all efficient estimators, they are therefore sensitive against misspecification of the model.

Specifically, consider observations $Y_1, \ldots, Y_n$ from a causal linear process

$$Y_t = X_t + \sum_{s=1}^{\infty} \delta_s X_{t-s}, \qquad t \in \mathbb{Z},$$

with independent and identically distributed innovations $X_t, t \in \mathbb{Z}$, with mean 0 and finite variance. A simple estimator of a linear functional $E[h(Y_0)]$ of the stationary distribution is the empirical estimator $\frac{1}{n}\sum_{j=1}^{n} h(Y_j)$. It does not use the fact that the process is linear and centered. We shall show how to construct better estimators if the process is invertible,

$$X_t = Y_t + \sum_{s=1}^{\infty} \gamma_s Y_{t-s}, \qquad t \in \mathbb{Z}.$$

The idea is to express the functional $E[h(Y_0)]$ as $E[h(X_0 + \sum_{s=1}^{\infty} \delta_s X_{-s})]$ and to estimate it by a $U$-statistic of increasing order based on estimated innovations, taking into account the constraint that the innovations have mean 0. We do this for a situation often encountered in applications: the coefficients $\delta_1, \delta_2, \ldots$ and hence also $\gamma_1, \gamma_2, \ldots$ depend on an unknown Euclidean parameter $\vartheta$.

The construction of our estimator involves several steps. Let us illustrate them with the simplest example, a linear autoregressive model of order 1,

$$Y_t = \vartheta Y_{t-1} + X_t, \qquad t \in \mathbb{Z},$$

with $\vartheta$ belonging to the interval $(-1, 1)$. Our result is new, and nontrivial, even for this simple case. The model is a semiparametric model with one-dimensional parameter $\vartheta$ and infinite-dimensional parameter $P$, the distribution of the innovations. The stationary distribution of this process thus depends on the pair $(\vartheta, P)$.

We want to estimate the linear functional $E[h(Y_0)]$ of the stationary distribution. The obvious estimator is again the empirical estimator $\frac{1}{n}\sum_{j=1}^{n} h(Y_j)$. It is known that the empirical estimator is a least dispersed regular estimator in Markov chain models with completely unspecified transition distribution; see Penev (1991), Bickel (1993) and Greenwood and Wefelmeyer (1995). Here, however, we are dealing with a semiparametric submodel. Thus, we should be able to improve upon this estimator.

Before we describe our estimator, let us briefly describe a simple improvement of the empirical estimator, obtained by exploiting the fact that the innovations, and hence the observations, have mean 0. This is a linear constraint $E[Y_0] = 0$ on the stationary distribution. For any $c \in \mathbb{R}$ we obtain a new estimator for $E[h(Y_0)]$:

$$\frac{1}{n}\sum_{j=1}^{n}(h(Y_j) - cY_j).$$

For general Markov chain models, Müller, Schick and Wefelmeyer (2001b) determine the constant $c$ which minimizes the asymptotic variance of the new estimator. For our autoregressive model, this constant becomes particularly simple if $h$ is a polynomial. For example, for the stationary *variance* $E[Y_0^2]$, that is, $h(y) = y^2$, the optimal constant is

$$c = c_* = \frac{\mu_3}{(1 + \vartheta)\mu_2},$$

with $\mu_k = E[X_1^k]$. This optimal $c_*$ depends on $P$ and $\vartheta$ and must be estimated. We estimate $\vartheta$ by the least squares estimator $\hat{\vartheta}_* = \frac{1}{n}\sum_{j=1}^n Y_{j-1}Y_j / \frac{1}{n}\sum_{j=1}^n Y_{j-1}^2$, the innovations by $Y_j - \hat{\vartheta}_* Y_{j-1}$ and $\mu_k$ by its empirical estimator based on estimated innovations:

$$\hat{\mu}_k = \frac{1}{n}\sum_{j=1}^n (Y_j - \hat{\vartheta}_* Y_{j-1})^k. \tag{1.1}$$

The resulting estimator for $E[Y_0^2]$ is

$$\frac{1}{n}\sum_{j=1}^n \left(Y_j^2 - \frac{\hat{\mu}_3}{(1+\hat{\vartheta}_*)\hat{\mu}_2} Y_j\right).$$

This simple improvement of the empirical estimator does not use the autoregressive structure of the chain. As mentioned above, this structure is exploited by a $U$-statistic of increasing order. Improving the empirical estimator then involves three steps. In the first step, we assume $\vartheta$ as known and exploit the structural relation $Y_t = \vartheta Y_{t-1} + X_t$. In the second step, we use the information that the innovation distribution has mean 0. The last step consists of replacing $\vartheta$ by an estimator.

The key step is the first one: we represent the observations as an infinite series of the innovations:

$$Y_t = \sum_{s=0}^\infty \vartheta^s X_{t-s}, \qquad t \in \mathbb{Z}.$$

Suppose first that the parameter $\vartheta$ is known. Then we can calculate the innovations $X_t = Y_t - \vartheta Y_{t-1}$, $t = 1, \ldots, n$, from the observations. Since $Y_0$ has the same distribution as $S = \sum_{s=1}^\infty \vartheta^{s-1} X_s$, the problem is now reduced to estimating the functional

$$E[h(Y_0)] = E[h(S)]$$
$$= E\left[h\left(\sum_{s=1}^\infty \vartheta^{s-1} X_s\right)\right]$$

from i.i.d. observations $X_1, \ldots, X_n$. This expectation is approximated by $E[h(S^{(m)})]$ with $S^{(m)} = \sum_{s=1}^m \vartheta^{s-1} X_s$ if $m$ increases with $n$. This suggests using the following variant of a $U$-statistic as an estimator for $E[h(S^{(m)})]$. Form the sums

$$S_i(\vartheta) = \sum_{s=1}^m \vartheta^{s-1} X_{i(s)} = \sum_{s=1}^m \vartheta^{s-1}(Y_{i(s)} - \vartheta Y_{i(s)-1})$$

for injective functions $i$ from $\{1, \ldots, m\}$ into $\{1, \ldots, n\}$. These sums are distributed as $S^{(m)}$. Hence we estimate $E[h(Y_0)]$ by an average over these

sums, the $U$-statistic

$$\hat{\kappa}(\vartheta) = \frac{(n-m)!}{n!} \sum_{i \in \Phi} h(S_i(\vartheta)),$$

where $\Phi$ denotes the set of all injective functions from $\{1, \ldots, m\}$ into $\{1, \ldots, n\}$.

We can show, via Hoeffding decomposition, that if $m = m(n)$ increases with $n$ at an appropriate rate, then the $U$-statistic $\hat{\kappa}(\vartheta)$ is *asymptotically linear*,

$$\hat{\kappa}(\vartheta) = E[h(Y_0)] + \frac{1}{n} \sum_{j=1}^{n} h_*(X_j) + o_p(n^{-1/2}),$$

with *influence function* $h_* = \sum_{s=1}^{\infty} h_s$, where $h_s(x) = E[h(S)|X_s = x] - E[h(S)]$. For fixed $m$, the $U$-statistic $\hat{\kappa}(\vartheta)$ is a least dispersed regular estimator of $E[h(S^{(m)})] = E[h(\sum_{s=1}^{m} \vartheta^{s-1} X_s)]$ if nothing is known about the distribution of the $X_j$. See Levit (1974), or argue via the asymptotic equivalence of the $U$-statistic and the von Mises statistic and efficiency of the empirical distribution function [Beran (1977)]. Optimality is preserved if we let $m$ tend to $\infty$ at the appropriate rate. For $U$-statistics of increasing order, see also Shieh (1994) and Heilig and Nolan (2001).

In Section 2 we prove these results for functionals of the more general form $E[h(\sum_{s=1}^{\infty} \beta_s X_s)]$ with summable coefficients $\beta_1, \beta_2, \ldots$. The results are of independent interest. For simplicity, we do not prove them under minimal assumptions on the function $h$. In our applications to linear time series in Sections 4 and 5, we shall need stronger assumptions anyway. The assumptions are general enough to cover moments and absolute moments and other smooth functions.

Now we turn to the second step of the construction of our estimator, exploiting the fact that $X_t$ has mean 0. This is a linear constraint of the form $E[Y_1 - \vartheta Y_0] = E[X_1] = 0$. The simple improvement of the empirical estimator $\frac{1}{n} \sum_{j=1}^{n} h(Y_j)$, described above, has used the linear constraint $E[Y_1] = 0$ on observations from a Markov chain. Here we use the constraint $E[X_1] = 0$ on the observed innovations, which are i.i.d. This simplifies improving our estimator $\hat{\kappa}(\vartheta)$. Similarly, as above, we form, for any $a \in \mathbb{R}$, the estimator

$$\hat{\kappa}(\vartheta, a) = \hat{\kappa}(\vartheta) - a \frac{1}{n} \sum_{j=1}^{n} (Y_j - \vartheta Y_{j-1}),$$

which has influence function $x \mapsto h_*(x) - ax$. It is easy to check that the choice

$$a = a_* = \frac{E[X_1 h_*(X_1)]}{E[X_1^2]}$$

yields an estimator with smallest asymptotic variance in this class of estimators. The optimal $a_*$ stems from projection on $[X_1]$. It depends on $P$ and must be replaced by an estimator. A consistent estimator is

$$\hat{a}_*(\vartheta) = \frac{\sum_{j=1}^{n}(Y_j - \vartheta Y_{j-1})\sum_{s=1}^{m} H_{s,j}(\vartheta)}{\sum_{j=1}^{n}(Y_j - \vartheta Y_{j-1})^2},$$

where

$$H_{s,j}(\vartheta) = \frac{(n-m)!}{(n-1)!} \sum_{i \in \Phi, i(s)=j} h(S_i(\vartheta)), \qquad s=1,\ldots,m, j=1,\ldots,n.$$

This leads us to the estimator

$$\hat{\kappa}(\vartheta, \hat{a}_*(\vartheta)) = \hat{\kappa}(\vartheta) - \hat{a}_*(\vartheta)\frac{1}{n}\sum_{j=1}^{n}(Y_j - \vartheta Y_{j-1}).$$

We show that this is a least dispersed regular estimator of $E[h(Y_0)]$ in the submodel with *known* parameter $\vartheta$. For a related efficiency result in such i.i.d. models with linear constraints, but for simpler functionals, see Levit (1975). In Section 3 we generalize these results to functionals of the form $E[h(\sum_{s=1}^{\infty} \beta_s X_s)]$.

The third and last step of the construction of our estimator consists of replacing $\vartheta$ by an estimator $\hat{\vartheta}$, leading to the *substitution estimator* $\hat{\kappa}(\hat{\vartheta}, \hat{a}_*(\hat{\vartheta}))$. It then follows from the *substitution principle* that the substitution estimator is efficient for $E[h(Y_0)] = E[h(\sum_{s=1}^{\infty} \vartheta^{s-1} X_s)]$ if $\hat{\vartheta}$ is efficient for $\vartheta$. Conditions for this principle to hold were first formulated by Klaassen and Putter (2001) in models with independent and identically distributed observations, and generalized to Markov chain models by Müller, Schick and Wefelmeyer (2001a).

In Section 4, rather than checking the conditions for the substitution principle, we calculate directly the influence function of the substitution estimator for functionals $E[h(\sum_{s=1}^{\infty} \alpha_s(\vartheta) X_s)]$ from observations which approximate $X_1,\ldots,X_n$. In Section 5 we apply the results of Sections 2–4 to estimate stationary expectations $E[h(Y_0)]$ from observations of causal invertible linear processes. Efficiency of our estimator follows from Schick and Wefelmeyer (2002a) who characterize efficient estimators for arbitrary differentiable functionals in such time series models.

In Section 6 we compare the asymptotic variances of the empirical estimator, the improved empirical estimator and our estimator for the stationary variance in AR(1) models. In this situation the asymptotic variances of the estimators can be calculated explicitly. For innovation distributions far from normal the variance decrease can be considerable.

**2. Estimating the distribution of an infinite series.** Let $X_1, X_2, \ldots$ be independent and identically distributed random variables with

$$E[|X_1|^{2p}] < \infty \tag{2.1}$$

for some $p \geq 1$ and with unknown common distribution $P$. Let $\beta_1, \beta_2, \ldots$ be known real numbers such that

$$\sum_{r=1}^{\infty} |\beta_r| < \infty. \tag{2.2}$$

Then the series

$$S = \sum_{r=1}^{\infty} \beta_r X_r$$

converges almost surely and in $L_{2p}$. Let $h$ be a function from $\mathbb{R}$ to $\mathbb{R}$ such that

$$|h(x)| \leq C_1(1 + |x|^p), \quad x \in \mathbb{R}, \tag{2.3}$$

$$|h(x+y) - h(x)| \leq C_2(1 + |x|^p)(|y| + |y|^p), \quad x, y \in \mathbb{R}, \tag{2.4}$$

for some finite constants $C_1$ and $C_2$. Then the expectation $E[h(S)]$ is well defined. Examples of functions $h$ that satisfy (2.3) and (2.4) are polynomials in $x$ or $|x|$ of degree at most $p$ and Lipschitz continuous functions.

We are interested in estimating $E[h(S)]$ from the observations $X_1, \ldots, X_n$. Let us introduce our estimator. It follows from (2.1)–(2.4) that the infinite sum $S$ is well approximated by the finite sum $S^{(m)} = \sum_{r=1}^{m} \beta_r X_r$ for moderately large $m$. Indeed, the Minkowski inequality yields that

$$E\left[\left|\sum_{j=a}^{b} \beta_j X_j\right|^q\right] \leq E[|X_1|^q]\left(\sum_{j=a}^{b} |\beta_j|\right)^q, \quad 1 \leq a \leq b, 1 \leq q \leq 2p. \tag{2.5}$$

In view of (2.4) and the independence of $S - S^{(m)}$ and $S^{(m)}$,

$$E[|h(S) - h(S^{(m)})|^2]$$
$$\leq C_2^2 E[(1 + |S^{(m)}|^p)^2]\Big(E[(|S - S^{(m)}| + |S - S^{(m)}|^p)^2]\Big).$$

It is now easy to see that there exists a constant $K$ such that

$$E[|h(S) - h(S^{(m)})|^2] \leq K^2 \left(\sum_{r=m+1}^{\infty} |\beta_r|\right)^2 \tag{2.6}$$

and hence

$$|E[h(S)] - E[h(S^{(m)})]| \leq K \sum_{r=m+1}^{\infty} |\beta_r|. \tag{2.7}$$

Actually, the constant $K$ can be chosen to be

$$K = 2C_2 \left(1 + \sum_{r=1}^{\infty} |\beta_r|\right)^{2p-1} (1 + E[X_1^2] + E[|X_1|^{2p}]).$$

Recall that $\Phi$ denotes the set of all injective functions from $\{1,\ldots,m\}$ to $\{1,\ldots,n\}$. The random variables

$$S_i = \sum_{r=1}^{m} \beta_r X_{i(r)}, \qquad i \in \Phi,$$

have the same distribution as $S^{(m)}$. Hence an unbiased estimator of $E[h(S^{(m)})]$ is given by

$$\tilde{\kappa} = \frac{(n-m)!}{n!} \sum_{i \in \Phi} h(S_i).$$

The estimator can be written as a $U$-statistic,

$$\binom{n}{m}^{-1} \sum_{1 \leq i(1) < \cdots < i(m) \leq n} k_m(X_{i(1)}, \ldots, X_{i(m)}),$$

with symmetric kernel $k_m$ defined by

$$k_m(x_1, \ldots, x_m) = \frac{1}{m!} \sum_{i \in \Pi} h(\beta_1 x_{i(1)} + \cdots + \beta_m x_{i(m)}), \qquad x_1, \ldots, x_m \in \mathbb{R},$$

with $\Pi$ the set of permutations of $\{1,\ldots,m\}$. Using standard $U$-statistic techniques [see Serfling (1980), page 178, Lemma A and page 184, Lemma B], we obtain

$$\tilde{\kappa} = \kappa_m + \frac{1}{n} \sum_{j=1}^{n} m k_{m,1}(X_j) + R,$$

where

$$\kappa_m = E[k_m(X_1, \ldots, X_m)] = E[h(S^{(m)})],$$
$$k_{m,1}(x) = E[k_m(x, X_2, \ldots, X_m)] - \kappa_m, \qquad x \in \mathbb{R},$$

and the remainder satisfies

$$E[R^2] \leq \sum_{r=2}^{m} \binom{m}{r}^2 \binom{n}{r}^{-1} E[k_m^2(X_1, \ldots, X_m)].$$

It is easy to check that $E[k_m^2(X_1,\ldots,X_m)] \leq E[h^2(S^{(m)})]$. Using $m!/(m-r)! \leq m^r$ and $n!/(n-r)! \geq (n-r)^r$, we obtain, for $n - m > m^2$,

$$E[R^2] \leq E[h^2(S^{(m)})] \sum_{r=2}^{m} \frac{1}{r!} \left(\frac{m^2}{n-m}\right)^r$$
$$\leq E[h^2(S^{(m)})] \left(\frac{m^2}{n-m}\right)^2.$$

Note also that

$$mk_{m,1}(x) = \sum_{r=1}^{m} \Big(E[h(S^{(m)})|X_r = x] - E[h(S^{(m)})]\Big), \qquad x \in \mathbb{R}.$$

Now let

$$h_r(x) = E[h(S)|X_r = x] - E[h(S)], \qquad x \in \mathbb{R}, r = 1, 2, \ldots.$$

With the help of (2.4) and the Cauchy–Schwarz inequality, we verify that

$$\int h_r^2 \, dP \leq 4E[(h(S) - h(S - \beta_r X_r))^2]$$

$$\leq 4C_2^2 E[(1 + |S - \beta_r X_r|^p)^2] E[(|\beta_r X_r| + |\beta_r X_r|^p)^2].$$

This and the Minkowski inequality show that there exists a constant $C$ such that, for all sufficiently large $m$ and $k$, $m < k$,

$$\int \Bigg(\sum_{r=m+1}^{k} h_r\Bigg)^2 dP \leq C \Bigg(\sum_{r=m+1}^{k} |\beta_r|\Bigg)^2.$$

Thus the series $h_* = \sum_{r=1}^{\infty} h_r$ is well defined in $L_2(P)$ and is the $L_2(P)$-limit of $\sum_{r=1}^{m} h_r$:

$$(2.8) \qquad \int \Bigg(h_* - \sum_{r=1}^{m} h_r\Bigg)^2 dP \to 0 \qquad \text{as } m \to \infty.$$

It follows from the Cauchy–Schwarz inequality and (2.6) that

$$(2.9) \qquad \int \Bigg(mk_{m,1} - \sum_{r=1}^{m} h_r\Bigg)^2 dP \leq 4mE[|h(S) - h(S^{(m)})|^2]$$

$$\leq 4mK^2 \Bigg(\sum_{r=m+1}^{\infty} |\beta_r|\Bigg)^2$$

for large $m$. We arrive at the following result.

THEOREM 2.1. *Suppose we can choose $m = m(n)$ such that*

$$(2.10) \qquad m^4/n \to 0 \quad \text{and} \quad n^{1/2} \sum_{r=m+1}^{\infty} |\beta_r| \to 0.$$

*Let $h$ satisfy* (2.3) *and* (2.4). *Then the estimator*

$$\tilde{\kappa} = \frac{(n-m)!}{n!} \sum_{i \in \Phi} h\Bigg(\sum_{r=1}^{m} \beta_r X_{i(r)}\Bigg)$$

is asymptotically linear for $E[h(S)]$ with influence function $h_* = \sum_{r=1}^{\infty} h_r$:

$$\tilde{\kappa} = E[h(S)] + \frac{1}{n}\sum_{j=1}^{n} h_*(X_j) + o_p(n^{-1/2}).$$

In particular, $\tilde{\kappa}$ is asymptotically normal with variance $\int h_*^2 \, dP$.

We have phrased this and the following theorems about estimators as asymptotic linearity results. The reason is that asymptotic linearity is useful for obtaining other, more familiar results about estimators: they are then seen to be asymptotically normal, their asymptotic variances are easily calculated and we can check whether they are regular and whether they are efficient in the sense of being least dispersed among regular estimators.

REMARK 2.1. Let us briefly discuss the choice of $m$ in two special cases:

1. Suppose that the coefficients $\beta_1, \beta_2, \ldots$ decay exponentially, say

$$|\beta_j| \leq C\vartheta^j, \qquad j = 1, 2, \ldots,$$

   for a finite constant $C$ and a positive number $\vartheta$, $\vartheta < 1$. Then the requirement (2.10) is satisfied if $m^4/n \to 0$ and $n^{1/2}\vartheta^m \to 0$. The latter holds if $\log(n)/m \to \infty$. If $\vartheta < e^{-1/2}$, it even holds for $m = \log(n)$.
2. Suppose $\beta_j = 0$ for $j > p$. Then we can take $m = p$. We should point out that in this case $h_* = h_1 + \cdots + h_p$ is a finite sum and (2.8) holds even though $m$ does not go to $\infty$. This is the classical result for fixed-degree $U$-statistics.

As it is very time consuming to calculate $\tilde{\kappa}$ for large $m$, it is advantageous to choose $m$ as small as possible.

REMARK 2.2. If the coefficients do not decay fast enough, we may not be able to satisfy (2.10). For example, if $\beta_j = j^{-1-a}$, $j = 1, 2, \ldots$, for some positive $a$, then $m$ needs to satisfy $m^4/n \to 0$ and $n/m^{2a} \to 0$. But this is only possible if $a > 2$.

Let us now show that $\tilde{\kappa}$ is efficient. For this it suffices to show that $E[h(S)]$ is differentiable at the true $P$ with canonical gradient equal to the influence function $h_*$ of our estimator $\tilde{\kappa}$. Since we will have to look at distributions near to, but different from, the true $P$, it will occasionally be convenient to express the dependence of expectations on the underlying distribution by writing $E_P$ for $E$. Note that $\kappa(P) = E_P[h(S)]$ defines a functional on the set of all distributions with finite $2p$th moments. We introduce a local model at the true $P$ as follows. Let $L_*(P)$ denote the set of all measurable functions $g$ from $\mathbb{R}$ to $\mathbb{R}$ such that $\int g\,dP = 0$ and $\int g^2\,dP < \infty$. To each $g$ in $L_*(P)$ associate a sequence $g_n$ in $L_*(P)$ such that

(2.11) $$|g_n| \leq n^{1/8} \quad \text{and} \quad \int (g_n - g)^2\,dP \to 0.$$

A possible choice is $g_n = g\mathbf{1}[2|g| \leq n^{1/8}] - \int g\mathbf{1}[2|g| \leq n^{1/8}]\,dP$. Let $P_{n,g}$ denote the distribution with $P$-density $1 + n^{-1/2}g_n$. Since $0 \leq 1 + n^{-1/2}g_n$ and $\int (1 + n^{-1/2}g_n)\,dP = 1$, the function $1 + n^{-1/2}g_n$ is indeed a probability density.

THEOREM 2.2. *Suppose we can choose $m = m(n)$ such that (2.10) holds. Let $h$ satisfy (2.3) and (2.4). Then the functional $\kappa(P) = E_P[h(S)]$ is differentiable at $P$ with gradient $h_* = \sum_{r=1}^{\infty} h_r$:*

$$n^{1/2}(\kappa(P_{n,g}) - \kappa(P)) \to \int h_* g\,dP.$$

PROOF. Let $m = m(n)$ satisfy (2.10). Let $G_{n,0} = 1$ and

$$G_{n,k} = \prod_{r=1}^{k}(1 + n^{-1/2}g_n(X_r)), \qquad k = 1, 2, \ldots.$$

Since

$$n^{1/2}(G_{n,k} - 1) = \sum_{r=1}^{k} G_{n,r-1}g_n(X_r) = \sum_{r=1}^{k} g_n(X_r) + \sum_{r=2}^{k} g_n(X_r)(G_{n,r-1} - 1)$$

and

$$E[(G_{n,k} - 1)^2] = E[G_{n,k}^2] - 1 = (1 + n^{-1}E[g_n^2(X_1)])^k - 1$$
$$\leq \frac{k}{n}E[g_n^2(X_1)](1 + n^{-1}E[g_n^2(X_1)])^{k-1},$$

we get by an application of the Cauchy–Schwarz inequality and the independence of $X_r$ and $G_{n,r-1}$ that

$$\left| n^{1/2}E[h(S)(G_{n,m} - 1)] - E\left[h(S)\sum_{r=1}^{m} g_n(X_r)\right] \right|$$
$$\leq \sum_{r=2}^{m}(E[h^2(S)]E[g_n^2(X_r)]E[(G_{n,r-1} - 1)^2])^{1/2} \to 0.$$

Since $\int g_n\,dP = 0$, we find that

$$E[h(S)g_n(X_r)] = E[(E[h(S)|X_r] - E[h(S)])g_n(X_r)]$$
$$= \int g_n h_r\,dP.$$

Thus, in view of (2.8) and (2.11),

$$E\left[h(S)\sum_{r=1}^{m}g_n(X_r)\right] = \int g_n \sum_{r=1}^{m} h_r\, dP \to \int g h_*\, dP.$$

This shows that

(2.12) $$n^{1/2}E[h(S)(G_{n,m}-1)] \to \int g h_*\, dP.$$

Note that $E_{P_{n,g}}[h(S^{(m)})] = E_P[h(S^{(m)})G_{n,m}]$, so that $\kappa(P_{n,g}) - \kappa(P)$ equals

$$E_{P_{n,g}}[h(S) - h(S^{(m)})] + E[(h(S^{(m)}) - h(S))G_{n,m}] + E[h(S)(G_{n,m}-1)].$$

The desired result now follows from (2.12) and (2.10) because

$$n^{1/2}|E_{P_{n,g}}[h(S) - h(S^{(m)})]| = O\left(n^{1/2}\sum_{r=m+1}^{\infty}|\beta_r|\right)$$

by the same argument that yields (2.7), and

$$n^{1/2}|E[(h(S) - h(S^{(m)}))G_{n,m}]| = O\left(n^{1/2}\sum_{r=m+1}^{\infty}|\beta_r|\right)$$

by (2.6) and $E[G_{n,m}^2] \to 1$. □

Theorems 2.1 and 2.2 imply that $\tilde{\kappa}$ is least dispersed among regular estimators of $E_P[h(S)]$ if nothing is known about $P$. For an appropriate version of the convolution theorem, see Bickel, Klaassen, Ritov and Wellner [(1998), page 63, Theorem 2, and page 65, Proposition 1].

**3. Estimation with constraints.** In the setting of Section 2, we can find better estimators for $E[h(S)] = E[h(\sum_{r=1}^{\infty}\beta_r X_r)]$ if additional information about the distribution $P$ is available. Suppose we know that

(3.1) $$\int \psi\, dP = 0$$

for some measurable function $\psi$ from $\mathbb{R}$ to $\mathbb{R}$ such that $\int \psi^2\, dP$ is finite and positive. An important case is the choice $\psi(x) = x$. This just means that $P$ has mean 0.

Under the constraint (3.1) we can consider the estimator

$$\tilde{\tilde{\kappa}}(a) = \tilde{\kappa} - a\frac{1}{n}\sum_{j=1}^{n}\psi(X_j)$$

for real $a$ and verify that it has influence function $h_* - a\psi$ if $m = m(n)$ satisfies (2.10):

$$\tilde{\tilde{\kappa}}(a) = E[h(S)] + \frac{1}{n}\sum_{j=1}^{n}(h_*(X_j) - a\psi(X_j)) + o_p(n^{-1/2}).$$

Its asymptotic variance is minimized for the choice

$$a = a_* = \frac{\int h_*\psi\,dP}{\int \psi^2\,dP},$$

which is the coefficient of the projection of $h_*$ onto $\psi$. Let us now construct an estimator of $a_*$ that is consistent if $m = m(n)$ satisfies (2.10). Our candidate is

$$\hat{a}_* = \frac{\sum_{j=1}^{n}\psi(X_j)\sum_{r=1}^{m}H_{r,j}}{\sum_{j=1}^{n}\psi^2(X_j)},$$

where

$$H_{r,j} = \frac{(n-m)!}{(n-1)!}\sum_{i\in\Phi, i(r)=j}h(S_i), \qquad r=1,\ldots,m, j=1,\ldots,n.$$

Recall that $S_i = \sum_{r=1}^{m}\beta_r X_{i(r)}$ for $i \in \Phi$. In view of the law of large numbers, we need only show that

(3.2) $$\frac{1}{n}\sum_{j=1}^{n}\psi(X_j)\sum_{r=1}^{m}H_{r,j} = \frac{1}{n}\sum_{j=1}^{n}h_*(X_j)\psi(X_j) + o_p(1).$$

Given $X_1$, the random variable $H_{r,1}$ is a $U$-statistic (of degree $m-1$ in the variables $X_2,\ldots,X_n$). Thus we have, for $r = 1,\ldots,m$ and $n - m \geq (m-1)^2$,

$$E[(H_{r,1} - E[H_{r,1}|X_1])^2]$$
$$\leq E[h^2(S^{(m)})]\sum_{k=1}^{m-1}\binom{m-1}{k}^2\binom{n-1}{k}^{-1}$$
$$\leq E[h^2(S^{(m)})]\frac{2(m-1)^2}{n-m}.$$

From this and the Cauchy–Schwarz inequality, we get

$$E\left[\frac{1}{n}\sum_{j=1}^{n}\left(\sum_{r=1}^{m}(H_{r,j} - E[H_{r,j}|X_j])\right)^2\right]$$
$$\leq m\sum_{r=1}^{m}E[(H_{r,1} - E[H_{r,1}|X_1])^2]$$
$$= O(m^4(n-m)^{-1}).$$

Thus $m^4/n \to 0$ implies that

$$\frac{1}{n}\sum_{j=1}^{n}\left(\sum_{r=1}^{m}(H_{r,j} - E[H_{r,j}|X_j])\right)^2 = o_p(1). \tag{3.3}$$

From this and another application of the Cauchy–Schwarz inequality, we can now conclude that

$$\frac{1}{n}\sum_{j=1}^{n}\psi(X_j)\sum_{r=1}^{m}H_{r,j} = \frac{1}{n}\sum_{j=1}^{n}\psi(X_j)\sum_{r=1}^{m}E[H_{r,j}|X_j] + o_p(1).$$

It is easy to check that

$$\sum_{r=1}^{m}E[H_{r,j}|X_j] = m(\kappa_m + k_{m,1}(X_j)). \tag{3.4}$$

As $m\kappa_m = o(n^{1/2})$, we obtain from the central limit theorem that

$$\frac{1}{n}\sum_{j=1}^{n}\psi(X_j)m\kappa_m = o_p(1).$$

In view of this, (2.8) and (2.9), we can now conclude the desired (3.2). Let us summarize this in the following theorem.

THEOREM 3.1. *Suppose we can choose* $m = m(n)$ *such that* (2.10) *holds. Let h satisfy* (2.3) *and* (2.4). *Then the estimator*

$$\tilde{\tilde{\kappa}}(\hat{a}_*) = \tilde{\kappa} - \hat{a}_*\frac{1}{n}\sum_{j=1}^{n}\psi(X_j)$$

*is asymptotically linear for* $\kappa(P) = E_P[h(S)]$ *with influence function* $h_* - a_*\psi$:

$$\tilde{\tilde{\kappa}}(\hat{a}_*) = \kappa(P) + \frac{1}{n}\sum_{j=1}^{n}[h_*(X_j) - a_*\psi(X_j)] + o_p(n^{-1/2}).$$

*In particular,* $\tilde{\tilde{\kappa}}(\hat{a}_*)$ *is asymptotically normal with variance*

$$E[(h_*(X_1) - a_*\psi(X_1))^2] = \int h_*^2\,dP - \frac{(\int h_*\psi\,dP)^2}{\int \psi^2\,dP}.$$

It is straightforward to check that $h_* - a_*\psi$ is the efficient influence function for estimators of $E_P[h(S)]$ under the constraint $\int \psi\,dP = 0$; see Levit (1975). It follows from Theorem 3.1 that $\tilde{\tilde{\kappa}}(\hat{a}_*)$ is a least dispersed regular estimator of $E_P[h(S)]$ when $P$ is unknown except for $\int \psi\,dP = 0$; see again the convolution theorem in Bickel, Klaassen, Ritov and Wellner [(1998), pages 63 and 65].

**4. Estimated coefficients and perturbed observations.** Let $X_1, \ldots, X_n$ be i.i.d. random variables with distribution $P$ satisfying (2.1). We want to estimate the expectation $E[h(\sum_{r=1}^{\infty} \beta_r X_r)]$. In the applications to time series we have in mind, the coefficients $\beta_1 = \alpha_1(\vartheta_0), \beta_2 = \alpha_2(\vartheta_0), \ldots$ depend on an unknown parameter $\vartheta_0$, and the random variables $X_1, \ldots, X_n$ are the unobservable innovations of a time series. In this case, both the coefficients and the innovations must be estimated from the time series using estimators of $\vartheta_0$. This will be done in Section 5. In preparation, the present section considers general estimators $X_{n,1}(\hat{\vartheta}), \ldots, X_{n,n}(\hat{\vartheta})$ of $X_1, \ldots, X_n$. Theorem 4.1 shows asymptotic linearity of a $U$-statistic based on observations $X_{n,1}(\hat{\vartheta}), \ldots, X_{n,n}(\hat{\vartheta})$; Theorem 4.2 treats the case with constraint $\int \psi \, dP = 0$. As the underlying parameter space we take an open subset $\Theta$ of $\mathbb{R}^d$. We assume that $\alpha_1, \alpha_2, \ldots$ are continuously differentiable functions from $\Theta$ to $\mathbb{R}$ such that, for some $\eta > 0$,

$$(4.1) \qquad \sum_{r=1}^{\infty} |\alpha_r(\vartheta_0)| < \infty \quad \text{and} \quad \sum_{r=1}^{\infty} \sup_{\|\vartheta - \vartheta_0\| < \eta} \|\dot{\alpha}_r(\vartheta)\| < \infty,$$

where $\dot{\alpha}_r$ denotes the gradient of $\alpha_r$. Note that this implies that

$$(4.2) \qquad \sum_{r=1}^{\infty} \sup_{\|\vartheta - \vartheta_0\| < \eta} |\alpha_r(\vartheta)| < \infty$$

for the same $\eta$ as in (4.1). We consider random variables $X_{n,1}(\vartheta), \ldots, X_{n,n}(\vartheta)$ such that $X_{n,j}(\vartheta)$ approximates $X_j$ if $\vartheta$ is close to $\vartheta_0$: there are $d$-dimensional random vectors $\xi_1, \xi_2, \ldots$ such that

$$(4.3) \qquad \sup_{j \geq 1} E[\|\xi_j\|^2] < \infty,$$

$$(4.4) \qquad \max_{j \leq n} n^{-1/2} \|\xi_j\| = o_p(1),$$

$$(4.5) \qquad \sup_{\|t\| \leq T} \sum_{j=1}^{n} (X_{n,j}(\vartheta_0 + n^{-1/2} t) - X_j - n^{-1/2} t^\top \xi_j)^2 = o_p(1)$$

for all finite $T$.

REMARK 4.1. Conditions (4.3) and (4.4) are implied by uniform integrability of the variables $\|\xi_1\|^2, \|\xi_2\|^2, \ldots$. The former is obvious; the latter follows as

$$P\left(\max_{1 \leq j \leq n} n^{-1/2} \|\xi_j\| > \eta\right)$$

$$\leq \frac{1}{n\eta^2} \sum_{j=1}^{n} E[\|\xi_j\|^2 \mathbf{1}[\|\xi_j\| > n^{1/2} \eta]]$$

$$\leq \frac{1}{\eta^2} \max_{1 \leq j \leq n} E[\|\xi_j\|^2 \mathbf{1}[\|\xi_j\| > n^{1/2} \eta]], \qquad \eta > 0.$$

Thus, if the random vectors $\xi_1, \xi_2, \ldots$ are identically distributed, then (4.3) and (4.4) follow from $E[\|\xi_1\|^2] < \infty$. Sufficient conditions for (4.5) are the asymptotic differentiability of $X_{n,j}$ at $\vartheta_0$ in the sense that

$$(4.6) \quad \sup_{\|t\| \leq T} \sum_{j=1}^n (X_{n,j}(\vartheta_0 + n^{-1/2}t) - X_{n,j}(\vartheta_0) - n^{-1/2}t^\top \dot{X}_{n,j}(\vartheta_0))^2 = o_p(1)$$

for all finite $T$ together with

$$(4.7) \quad \frac{1}{n} \sum_{j=1}^n \|\dot{X}_{n,j}(\vartheta_0) - \xi_j\|^2 = o_p(1),$$

$$(4.8) \quad \sum_{j=1}^n (X_{n,j}(\vartheta_0) - X_j)^2 = o_p(1).$$

In applications to time series, $X_{n,j}(\vartheta_0)$ is a truncated series representation of innovations; see (5.5).

For $\vartheta \in \Theta$ and $i \in \Phi$, set now

$$S(\vartheta) = \sum_{r=1}^\infty \alpha_r(\vartheta) X_r,$$

$$S_i(\vartheta) = \sum_{r=1}^m \alpha_r(\vartheta) X_{i(r)},$$

$$S_{n,i}(\vartheta) = \sum_{r=1}^m \alpha_r(\vartheta) X_{n,i(r)}(\vartheta).$$

Set $S = S(\vartheta_0)$ and $S_i = S_i(\vartheta_0)$. These are the series in Section 2. Think of $S_{n,i}(\vartheta)$ as an approximation of $S_i(\vartheta)$. Next define

$$\hat{\kappa}(\vartheta) = \frac{(n-m)!}{n!} \sum_{i \in \Phi} h(S_{n,i}(\vartheta)), \qquad \vartheta \in \Theta.$$

Then $\hat{\kappa}(\vartheta_0)$ is an "estimator" of $E[h(S)]$ and defined as in Section 2, but now with $X_1, \ldots, X_n$ replaced by $X_{n,1}(\vartheta_0), \ldots, X_{n,n}(\vartheta_0)$. Let $\hat{\vartheta}$ be an estimator of $\vartheta_0$. In this section we calculate the influence function of $\hat{\kappa}(\hat{\vartheta})$. The result will be used in Section 5.

ASSUMPTION H. The function $h$ satisfies (2.3) and (2.4) and is absolutely continuous with an almost everywhere derivative $h'$ that is almost surely continuous with respect to the distribution of $S$ and satisfies the growth condition

$$|h'(x)| \leq C_3(1 + |x|)^q, \qquad x \in \mathbb{R},$$

for some constant $C_3$ and some $q \in [0, p]$.

Examples of functions $h$ that satisfy Assumption H are again polynomials in $x$ or $|x|$ of degree at most $p$ and Lipschitz continuous functions.

THEOREM 4.1. *Suppose assumptions* (4.1)–(4.5) *hold, $h$ satisfies Assumption* H *and we can choose $m = m(n)$ such that* (2.10) *holds with $\beta_r = \alpha_r(\vartheta_0)$. If $\hat{\vartheta}$ is $n^{1/2}$-consistent for $\vartheta_0$, then*

$$\hat{\kappa}(\hat{\vartheta}) = \tilde{\kappa} + \Delta_n^\top(\hat{\vartheta} - \vartheta_0) + o_p(n^{-1/2}), \tag{4.9}$$

*where*

$$\Delta_n = \frac{(n-m)!}{n!} \sum_{i \in \Phi} h'(S_i) D_i,$$

$$D_i = \sum_{r=1}^m [\dot{\alpha}_r(\vartheta_0) X_{i(r)} + \alpha_r(\vartheta_0) \xi_{i(r)}], \qquad i \in \Phi.$$

PROOF. For $i \in \Phi$ set

$$D_{n,i} = S_{n,i}(\hat{\vartheta}) - S_i = \sum_{r=1}^m [\alpha_r(\hat{\vartheta}) X_{n,i(r)}(\hat{\vartheta}) - \alpha_r(\vartheta_0) X_{i(r)}].$$

Since $h$ is absolutely continuous, we see that

$$\hat{\kappa}(\hat{\vartheta}) - \tilde{\kappa} = \frac{(n-m)!}{n!} \sum_{i \in \Phi} D_{n,i} \int_0^1 h'(S_i + z D_{n,i}) \, dz.$$

The desired result can now be written as

$$\frac{(n-m)!}{n!} \sum_{i \in \Phi} \left( D_{n,i} \int_0^1 h'(S_i + z D_{n,i}) \, dz - D_i^\top (\hat{\vartheta} - \vartheta_0) h'(S_i) \right) = o_p(n^{-1/2}).$$

But this is a consequence of the following statements:

$$\frac{(n-m)!}{n!} \sum_{i \in \Phi} (h'(S_i))^2 = O_p(1), \tag{4.10}$$

$$\frac{(n-m)!}{n!} \sum_{i \in \Phi} \|D_i\|^2 = O_p(1), \tag{4.11}$$

$$\frac{(n-m)!}{n!} \sum_{i \in \Phi} (D_{n,i} - D_i^\top (\hat{\vartheta} - \vartheta_0))^2 = o_p(n^{-1}), \tag{4.12}$$

$$\frac{(n-m)!}{n!} \sum_{i \in \Phi} \int_0^1 (h'(S_i + z D_{n,i}) - h'(S_i))^2 \, dz = o_p(1). \tag{4.13}$$

Of course, (4.10) holds because its left-hand side has an expectation that converges to that of $E[h'(S)^2]$ by the properties of $h'$. Next, we have

$$E[\|D_i\|^2] \leq 2\left(\sum_{r=1}^{m} \|\dot{\alpha}_r(\vartheta_0)\|\right)^2 E[X_1^2]$$
$$+ 2\left(\sum_{r=1}^{m} |\alpha_r(\vartheta_0)|\right)^2 \max_{1\leq j\leq n} E[\|\xi_j\|^2] \quad (4.14)$$

by the following version of the Cauchy–Schwarz inequality:

$$\left(\sum_r a_r b_r\right)^2 \leq \sum_r |a_r| \sum_r |a_r| b_r^2.$$

Relation (4.11) follows from (4.14) and assumptions (4.1)–(4.3). To obtain relation (4.12) use the formula

$$\frac{(n-m)!}{n!} \sum_{i\in\Phi} \left(\sum_{r=1}^{m} |a_r b_{i(r)}|\right)^2 \leq \left(\sum_{r=1}^{m} |a_r|\right)^2 \frac{1}{n}\sum_{j=1}^{n} b_j^2$$

to bound the left-hand side of (4.12) by

$$3\left(\sum_{r=1}^{m} |\alpha_r(\hat{\vartheta}) - \alpha_r(\vartheta_0)|\right)^2 \frac{1}{n}\sum_{j=1}^{n} \|\xi_j\|^2 \|\hat{\vartheta}-\vartheta_0\|^2$$
$$+ 3\left(\sum_{r=1}^{m} |\alpha_r(\hat{\vartheta})|\right)^2 \frac{1}{n}\sum_{j=1}^{n}(X_{n,j}(\hat{\vartheta}) - X_j - \xi_j^\top(\hat{\vartheta}-\vartheta_0))^2$$
$$+ 3\left(\sum_{r=1}^{m} \int_0^1 \|\dot{\alpha}_r(\vartheta_0 + z(\hat{\vartheta}-\vartheta_0)) - \dot{\alpha}_r(\vartheta_0)\| dz\right)^2 \frac{1}{n}\sum_{j=1}^{n} |X_j|^2 \|\hat{\vartheta}-\vartheta_0\|^2.$$

The desired (4.12) is now immediate in view of (4.1)–(4.5) and the $n^{1/2}$-consistency of $\hat{\vartheta}$. Note that the $n^{1/2}$-consistency of $\hat{\vartheta}$, the continuity of $\dot{\alpha}_r$ and (4.1) yield

$$\sum_{r=1}^{m} \int_0^1 \|\dot{\alpha}_r(\vartheta_0 + z(\hat{\vartheta}-\vartheta_0)) - \dot{\alpha}_r(\vartheta_0)\| dz = o_p(1).$$

We also have

$$D_n = \max_{i\in\Phi} |D_{n,i}| = o_p(1). \quad (4.15)$$

This is a consequence of (4.12) and the fact that

$$\max_{i\in\Phi} n^{-1/2}\|D_i\| \leq \sum_{r=1}^{\infty} \|\dot{\alpha}_r(\vartheta_0)\| \max_{1\leq j\leq n} n^{-1/2}|X_j|$$

$$+ \sum_{r=1}^{\infty} |\alpha_r(\vartheta_0)| \max_{1 \le j \le n} n^{-1/2} \|\xi_j\|$$
$$= o_p(1).$$

Thus it suffices to prove (4.13) with $D_{n,i}$ replaced by $D_{n,i}^* = D_{n,i}\mathbf{1}[|D_{n,i}| \le 1]$. It follows from Assumption H that

$$Z_{n,i} = \int_0^1 (h'(S_i + zD_{n,i}^*) - h'(S_i))^2 \, dz \le 4C_3^2(2 + |S_i|)^{2p}, \qquad i \in \Phi.$$

Since $S_i$ has the same distribution as $S^{(m)}(\vartheta_0) = \sum_{r=1}^m \alpha_r(\vartheta_0)X_j$ and $S^{(m)}$ converges in $L_{2p}$ to $S$, we see that the random variables $\{Z_{n,i} : i \in \Phi, n \ge 1\}$ are uniformly integrable. Thus (4.13) will follow if we can show that, for every $L$,

(4.16) $$\frac{(n-m)!}{n!} \sum_{i \in \Phi} \int_0^1 L \wedge (h'(S_i + zD_{n,i}^*) - h'(S_i))^2 \, dz = o_p(1).$$

Fix $L$. Define a map $H$ from $\mathcal{Q}$, the set of all probability measures on the Borel $\sigma$-field of $\mathbb{R}^2$, into $[0, L]$ by

$$H(Q) = \int L \wedge (h'(x) - h'(y))^2 Q(dx, dy), \qquad Q \in \mathcal{Q}.$$

With the aid of this map, we can write the expected value of the left-hand side of (4.16) as

$$\frac{(n-m)!}{n!} \sum_{i \in \Phi} \int_0^1 H(Q_{n,i}^z) \, dz,$$

where $Q_{n,i}^z$ is the distribution of the bivariate random vector

$$Y_{n,i}^z = (S_i + zD_{n,i}^*, S_i)^\top.$$

Endow $\mathcal{Q}$ with the topology of weak convergence. This topology is generated by the Prohorov metric $\rho$. By the properties of $h'$, the map $H$ is bounded and continuous at $Q_0$, the distribution of $(S, S)^\top$. Note also that $H(Q_0) = 0$. Hence for $\varepsilon > 0$ there exists $\delta > 0$ such that $\rho(Q, Q_0) < \delta$ implies $|H(Q)| < \varepsilon$. It thus suffices to show that

(4.17) $$\sup\{\rho(Q_{n,i}^z, Q_0) : i \in \Phi, z \in [0,1]\} \to 0.$$

For this we use the following simple property of the Prohorov metric. If $X$ and $Y$ are two bivariate random vectors with distributions $Q$ and $R$, then $\rho(Q, R) \le \eta + P(\|X - Y\| \ge \eta)$ for each $\eta > 0$. Now let $Y_i = (S_i^*, S_i^*)^\top$ with

$$S_i^* = S_i + \sum_{r=1}^{\infty} \alpha_{m+r}(\vartheta_0) X_{n+r}.$$

Then $Y_i$ has distribution $Q_0$ and $\|Y_{n,i}^z - Y_i\| \le \sqrt{2}|S_i - S_i^*| + |D_{n,i}^*|$ for all $z \in [0,1]$ and all $i \in \Phi$. The desired (4.17) is now immediate. $\square$

REMARK 4.2. For $i \in \Phi$ set

$$\dot{S}_i = \sum_{r=1}^{m} \dot{\alpha}_r(\vartheta_0) X_{i(r)} \quad \text{and} \quad T_i = \sum_{r=1}^{m} \alpha_r(\vartheta_0) \xi_{i(r)},$$

so that $D_i = \dot{S}_i + T_i$. Under the assumptions of Theorem 4.1, one can show that

(4.18) $$\frac{(n-m)!}{n!} \sum_{i \in \Phi} h'(S_i)\dot{S}_i = \mu + o_p(1),$$

where $\mu = E[h'(S)\dot{S}]$ with $\dot{S} = \sum_{r=1}^{\infty} \dot{\alpha}_r(\vartheta_0) X_r$.

One also expects that under mild additional assumptions,

(4.19) $$\frac{(n-m)!}{n!} \sum_{i \in \Phi} h'(S_i)T_i = \nu + o_p(1)$$

for some vector $\nu \in \mathbb{R}^d$. Then (4.9) simplifies to

$$\hat{\kappa}(\hat{\vartheta}) = \tilde{\kappa} + (\mu + \nu)^\top (\hat{\vartheta} - \vartheta_0) + o_p(n^{-1/2}).$$

In the following lemma we formulate a set of sufficient conditions for (4.19) that is useful for the applications we have in mind.

LEMMA 4.1. *Suppose Assumption* H *holds, $m^4/n \to 0$, the random vectors $\xi_1, \xi_2, \ldots$ are stationary with $E[\|\xi_1\|^2] < \infty$ and*

(4.20) $$\sup_{r>s} E[\|\xi_r - E[\xi_r | X_{r-s}, \ldots, X_{r-1}]\|^2] \to 0 \quad \text{as } s \to \infty.$$

*Then* (4.19) *holds with*

$$\nu = E[h'(S)] \sum_{r=1}^{\infty} \alpha_r(\vartheta_0) E[\xi_1].$$

PROOF. Without loss of generality, we may assume that $d = 1$. Let $s$ denote the integer part of $1 + \log(n)$. Let $\Phi_s$ denote the set of all $i$ in $\Phi$ such that $i(q) > s$ and $|i(q) - i(r)| > s$ for all $q, r = 1, \ldots, m$ and $q \neq r$. Set $\xi_{r,s} = E[\xi_r | X_{r-s}, \ldots, X_{r-1}]$ for $r > s$ and

$$T_{i,s} = \sum_{r=1}^{m} \alpha_r(\vartheta_0) \xi_{i(r),s}, \quad i \in \Phi.$$

Since $m^4/n \to 0$, we have that $n^m/(n-ms)^m \to 1$. This shows that the cardinality of $\Phi_s$ is of the same order as that of $\Phi$. Hence the cardinality of the complement $\Phi \setminus \Phi_s$ of $\Phi_s$ with respect to $\Phi$ is of order $o(n!/(n-m)!)$.

We now use this and (4.20) to show that the left-hand side of (4.19) differs from

$$D = \frac{(n-m)!}{n!} \sum_{i \in \Phi_s} h'(S_i) T_{i,s}$$

by a term of order $o_p(1)$. Indeed, the expected value of the absolute value of this term is bounded by

$$\frac{(n-m)!}{n!} \left( \sum_{i \in \Phi \setminus \Phi_s} E[|h'(S_i) T_i|] + \sum_{i \in \Phi_s} E[|h'(S_i)(T_{i,s} - T_i)|] \right).$$

Now use the fact that the expected values $E[h'(S_i)^2]$ and $E[T_i^2]$ are uniformly bounded and that $E[(T_{i,s} - T_i)^2]^{1/2} \leq \sum_{j=1}^{\infty} |\alpha_j(\vartheta_0)| \sup_{r>s}(E[\|\xi_r - \xi_{r,s}\|^2])^{1/2}$, to conclude that this bound tends to 0.

It is easy to check that two summands $h'(S_i)T_{i,s}$ and $h'(S_j)T_{j,s}$ of $D$ are independent if their indices $i$ and $j$ satisfy $|i(r) - j(r)| > s$ for all $r = 1, \ldots, m$. This shows that the variance of $D$ goes to 0, so that $D = E[D] + o_p(1)$. Since $S_i$ and $T_{i,s}$ are independent for $i \in \Phi_s$, and $S_i$ has the same distribution as $S^{(m)} = \sum_{r=1}^{m} \alpha_r(\vartheta_0) X_r$, we have

$$E[D] = E[h'(S^{(m)})] \frac{(n-m)!}{n!} \sum_{i \in \Phi_s} E[T_{i,s}].$$

The properties of $h'$ imply $E[h'(S^{(m)})] \to E[h'(S)]$. From (4.20) and (4.1), we get

$$\sup_{i \in \Phi_s} \left| E[T_{i,s}] - \sum_{r=1}^{\infty} \alpha_r(\vartheta_0) E[\xi_1] \right| \to 0.$$

We can now conclude that $E[D] \to \nu$. This completes the proof. $\square$

Let us now turn to the constrained setting of Section 3, with $\psi$ a function such that $\int \psi \, dP = 0$ and $\int \psi^2 \, dP$ finite and positive. For $\vartheta \in \Theta$ consider

$$\hat{\hat{\kappa}}(\vartheta, \hat{a}_*(\vartheta)) = \hat{\kappa}(\vartheta) - \hat{a}_*(\vartheta) \frac{1}{n} \sum_{j=1}^{n} \psi(X_{n,j}(\vartheta)),$$

where

$$\hat{a}_*(\vartheta) = \frac{\sum_{j=1}^{n} \psi(X_{n,j}(\vartheta)) \sum_{r=1}^{m} H_{r,j}(\vartheta)}{\sum_{j=1}^{n} \psi^2(X_{n,j}(\vartheta))},$$

$$H_{r,j}(\vartheta) = \frac{(n-m)!}{(n-1)!} \sum_{i \in \Phi, i(r)=j} h(S_i(\vartheta)), \qquad r = 1, \ldots, m, j = 1, \ldots, n.$$

We now write $a_*(\vartheta_0)$ for the $a_*$ of Section 3 to stress the dependence on the parameter.

THEOREM 4.2. *Suppose the assumptions of Theorem* 4.1 *hold. Suppose also that $\psi$ is Lipschitz with an almost everywhere derivative $\psi'$ that is continuous P-almost surely. If $\hat{\vartheta}$ is $n^{1/2}$-consistent for $\vartheta_0$, then*

$$\frac{1}{n}\sum_{j=1}^{n} \hat{a}_*(\hat{\vartheta})\psi(X_{n,j}(\hat{\vartheta}))$$

$$= a_*(\vartheta_0)\frac{1}{n}\sum_{j=1}^{n} \psi(X_j) + a_*(\vartheta_0)\Gamma_n^\top(\hat{\vartheta} - \vartheta_0) + o_p(n^{-1/2}),$$

*with*

$$\Gamma_n = \frac{1}{n}\sum_{j=1}^{n} \psi'(X_j)\xi_j.$$

*If, in addition, the random vectors $\xi_1, \xi_2, \ldots$ are stationary and satisfy* (4.20), *then*

(4.21) $$\Gamma_n = E[\psi'(X_1)]E[\xi_1] + o_p(1),$$

*and hence $\hat{\hat{\kappa}}_n^* = \hat{\hat{\kappa}}(\hat{\vartheta}, \hat{a}_*(\hat{\vartheta}))$ equals*

$$\hat{\hat{\kappa}}_n^* = \hat{\kappa}(\hat{\vartheta}) - a_*(\vartheta_0)\frac{1}{n}\sum_{j=1}^{n} \psi(X_j)$$

$$- a_*(\vartheta_0)E[\psi'(X_1)]E[\xi_1^\top](\hat{\vartheta} - \vartheta_0) + o_p(n^{-1/2}).$$

PROOF. Since (4.21) is easy, we prove only the first conclusion. It suffices to show that

(4.22) $$\hat{a}_*(\hat{\vartheta}) = a_*(\vartheta_0) + o_p(1),$$

(4.23) $$\frac{1}{n}\sum_{j=1}^{n}(\psi(X_{n,j}(\hat{\vartheta})) - \psi(X_j)) = \frac{1}{n}\sum_{j=1}^{n} \psi'(X_j)\xi_j^\top(\hat{\vartheta} - \vartheta_0) + o_p(n^{-1/2}).$$

The latter is a special case of Theorem 4.1 with $h$ replaced by $\psi$ and $\alpha_1(\vartheta) = 1$ and $\alpha_r(\vartheta) = 0$ for $r \geq 2$. As $\psi$ is Lipschitz, we obtain from (4.4), (4.5) and the $n^{1/2}$-consistency of $\hat{\vartheta}$ that

(4.24) $$\Psi_n = \max_{1 \leq j \leq n}|\psi(X_{n,j}(\hat{\vartheta})) - \psi(X_j)| = o_p(1).$$

In view of (3.2) and (4.24), the desired statement (4.22) will follow from

(4.25) $$\frac{1}{n}\sum_{j=1}^{n}(\psi(X_{n,j}(\hat{\vartheta})) - \psi(X_j))\sum_{r=1}^{m} H_{r,j}(\vartheta_0) = o_p(1),$$

(4.26) $$\frac{1}{n}\sum_{j=1}^{n}\psi(X_{n,j}(\hat{\vartheta}))\sum_{r=1}^{m}(H_{r,j}(\hat{\vartheta}) - H_{r,j}(\vartheta_0)) = o_p(1).$$

It follows from (2.9), (3.3) and (3.4) that

$$\frac{1}{n}\sum_{j=1}^{n}\left(m\kappa_m - \sum_{r=1}^{m} H_{r,j}(\vartheta_0)\right)^2 = O_p(1).$$

It follows from (4.23) that

$$\frac{1}{n}\sum_{j=1}^{n}(\psi(X_{n,j}(\hat{\vartheta})) - \psi(X_j))m\kappa_m = o_p(1).$$

Together with (4.24), these statements yield (4.25). Next bound the absolute value of the left-hand side of (4.26) by

$$\sum_{r=1}^{m}\frac{(n-m)!}{n!}\sum_{j=1}^{n}\sum_{i\in\Phi, i(r)=j}|\psi(X_{n,j}(\hat{\vartheta}))|C_3(1+|S_i|+|D_{n,i}|)^q|D_{n,i}|$$

$$\leq C_4 \sum_{r=1}^{m}\frac{(n-m)!}{n!}\sum_{j=1}^{n}\sum_{i\in\Phi, i(r)=j}(1+|X_j|+\Psi_n)(1+|S_i|+D_n)^q|D_{n,i}|,$$

where $D_{n,i}$ and $D_n$ are as in the proof of Theorem 4.1 and $C_4$ is a constant. An application of the Cauchy–Schwarz inequality now shows that the square of the left-hand side of (4.26) is bounded by $m^2 C^2 U_n V_n$, where

$$U_n = \frac{(n-m)!}{n!}\sum_{i\in\Phi}D_{n,i}^2,$$

$$V_n = \frac{1}{m}\sum_{r=1}^{m}\frac{(n-m)!}{n!}\sum_{j=1}^{n}\sum_{i\in\Phi, i(r)=j}(1+|X_j|+\Psi_n)^2(1+|S_i|+D_n)^{2q}.$$

It follows from (4.11) and (4.12) that $nU_n = O_p(1)$. It follows from (4.24), (4.15) and $q \leq p-1$ that $V_n = O_p(1)$. As $m^2/n \to 0$, we obtain the desired (4.26). □

**5. Application to semiparametric linear processes.** Now we apply Sections 2–4 to real-valued causal invertible processes $Y_t, t \in \mathbb{Z}$, with infinite-order moving average and autoregressive representations

$$(5.1) \qquad Y_t = X_t + \sum_{s=1}^{\infty}\delta_s(\vartheta)X_{t-s}, \qquad t \in \mathbb{Z},$$

$$(5.2) \qquad Y_t = X_t - \sum_{s=1}^{\infty}\gamma_s(\vartheta)Y_{t-s}, \qquad t \in \mathbb{Z},$$

where the innovations $\{X_t, t \in \mathbb{Z}\}$ are i.i.d. with distribution $P$ which has mean 0 and finite variance, and the parameter $\vartheta$ varies in an open subset $\Theta$ of $\mathbb{R}^d$. We assume that $\delta_1, \delta_2, \ldots$ and $\gamma_1, \gamma_2, \ldots$ are continuously differentiable functions from $\Theta$ into $\mathbb{R}$ with the following growth conditions at the true parameter $\vartheta = \vartheta_0$: for a finite constant $C$ and positive numbers $\eta$ and $a < 1$,

$$\sup_{\|\vartheta - \vartheta_0\| < \eta} [|\delta_r(\vartheta)| + \|\dot{\delta}_r(\vartheta)\|] \le Ca^r, \qquad r = 1, 2, \ldots, \tag{5.3}$$

$$\sup_{\|\vartheta - \vartheta_0\| < \eta} [|\gamma_r(\vartheta)| + \|\dot{\gamma}_r(\vartheta)\|] \le Ca^r, \qquad r = 1, 2, \ldots. \tag{5.4}$$

Here $\dot{\delta}_r$ is the gradient of $\delta_r$, and $\dot{\gamma}_r$ the gradient of $\gamma_r$.

EXAMPLE 5.1. For the AR(1) process $Y_t = X_t + \vartheta Y_{t-1}$, take $\Theta = (-1, 1)$ and set $\gamma_1(\vartheta) = -\vartheta$ and $\gamma_s(\vartheta) = 0$ for $s \ge 2$. The infinite-order moving average representation holds with $\delta_s(\vartheta) = \vartheta^s$.

EXAMPLE 5.2. For the MA(1) process $Y_t = X_t + \vartheta X_{t-1}$, take $\Theta = (-1, 1)$ and set $\delta_1(\vartheta) = \vartheta$ and $\delta_s(\vartheta) = 0$ for $s \ge 2$. The infinite-order autoregressive representation holds with $\gamma_s(\vartheta) = (-\vartheta)^s$.

EXAMPLE 5.3. For the ARMA(1,1) process $Y_t - \vartheta_1 Y_{t-1} = X_t - \vartheta_2 X_{t-1}$, take $\Theta = \{(\vartheta_1, \vartheta_2) : \vartheta_1, \vartheta_2 \in (-1, 1), \vartheta_1 \ne \vartheta_2\}$. The infinite-order moving average representation holds with $\delta_s(\vartheta) = (\vartheta_1 - \vartheta_2)\vartheta_1^{s-1}$, and the infinite-order autoregressive representation holds with $\gamma_s(\vartheta) = (\vartheta_2 - \vartheta_1)\vartheta_2^{s-1}$.

In the following, we will occasionally write $Y_t(\vartheta)$ for representation (5.1) of $Y_t$, and $E_P$ for expectation when $P$ is true. We want to estimate the functional

$$\kappa(\vartheta, P) = E_P[h(Y_1(\vartheta))]$$

from observations $Y_0, \ldots, Y_n$. Since the true innovation distribution $P$ has mean 0, we have the linear constraint $\int x\, P(dx) = 0$, that is, $E_P[\psi(X_1)] = 0$ for $\psi(x) = x$.

Note that if we observe only $Y_1, \ldots, Y_n$, we cannot estimate the first few innovations so well that (4.5) holds. However, (4.5) can be achieved if we also observe $Y_{-r(n)}, \ldots, Y_0$ for a properly chosen sequence $r(n)$ of integers. For example, $r(n) = p - 1$ works for AR($p$). In general, we must have Assumption 3 in Schick and Wefelmeyer (2002a), which under our assumption (5.4) holds with $r(n)$ proportional to $(\log n)^{1+\varepsilon}$ for some $\varepsilon > 0$. We will assume in this section that those additional observations are available. Otherwise, renumber the observations.

We apply Section 4 with $\alpha_r = \delta_{r-1}$, $r = 1, 2, \ldots$, where $\delta_0 = 1$, and take $X_{n,1}(\vartheta), \ldots, X_{n,n}(\vartheta)$ to be truncated versions of the representation (5.2) of the innovations $X_1, \ldots, X_n$ in terms of the observations:

$$X_{n,j}(\vartheta) = Y_j + \sum_{s=1}^{r(n)+j} \gamma_s(\vartheta) Y_{j-s}, \qquad j = 1, \ldots, n,\ \vartheta \in \Theta. \tag{5.5}$$

It is easy to see that assumption (5.3) implies assumptions (4.1) and (4.2). Let us now show that (5.3) and (5.4) imply (4.3)–(4.5) with

$$\xi_j = \sum_{s=1}^{\infty} \dot{\gamma}_s(\vartheta_0) Y_{j-s}, \qquad j = 1, 2, \dots.$$

As $\xi_1, \xi_2, \dots$ are stationary and square integrable by (5.4), we obtain (4.3) and (4.4) from Remark 4.1. To prove relation (4.5), we verify the sufficient conditions (4.6)–(4.8) with $\dot{X}_{n,j}(\vartheta) = \sum_{s=1}^{r(n)+j} \dot{\gamma}_s(\vartheta) Y_{j-s}$. Conditions (4.7) and (4.8) are easy consequences of the choice of $r(n)$ and assumption (5.4). We bound the expectation of the left-hand side of (4.6) by

$$\sum_{j=1}^{n} E \left( \sum_{s=1}^{r(n)+j} \sup_{\|t\| \leq T} |\gamma_s(\vartheta_0 + n^{-1/2} t) - \gamma_s(\vartheta_0) - n^{-1/2} t^{\top} \dot{\gamma}_s(\vartheta_0)| |Y_{j-s}| \right)^2$$

$$\leq E(Y_1^2) \left( \sum_{s=1}^{\infty} \sup_{\|t\| \leq T} |\gamma_s(\vartheta_0 + n^{-1/2} t) - \gamma_s(\vartheta_0) - n^{-1/2} t^{\top} \dot{\gamma}_s(\vartheta_0)| \right)^2.$$

We have used the Minkowski inequality here. Since $\gamma_1, \gamma_2, \dots$ are continuously differentiable, each term in the last series converges to 0 as $n$ tends to $\infty$. Hence the sequence of series converges to 0 since the dominated convergence theorem applies by (5.4). This proves (4.6) and completes the proof of (4.5). Finally, assumptions (5.3) and (5.4) imply relation (4.20).

Now set

$$S_{n,i}(\vartheta) = \sum_{r=1}^{m} \delta_{r-1}(\vartheta) X_{n,i(r)}(\vartheta),$$

$$H_{r,j}(\vartheta) = \frac{(n-m)!}{(n-1)!} \sum_{i \in \Phi, i(r)=j} h(S_{n,i}(\vartheta)),$$

$$\hat{a}_*(\vartheta) = \frac{\sum_{j=1}^{n} (Y_j - \vartheta Y_{j-1}) \sum_{r=1}^{m} H_{r,j}(\vartheta)}{\sum_{j=1}^{n} (Y_j - \vartheta Y_{j-1})^2},$$

$$\hat{\hat{\kappa}}(\vartheta, a) = \frac{(n-m)!}{n!} \sum_{i \in \Phi} h(S_{n,i}(\vartheta)) - a \frac{1}{n} \sum_{j=1}^{n} (Y_j - \vartheta Y_{j-1}).$$

Since the random vectors $\xi_1, \xi_2, \dots$ are stationary with $E_P[\xi_1] = 0$, Theorem 4.2 implies

$$\hat{\hat{\kappa}}(\hat{\vartheta}, \hat{a}_*(\hat{\vartheta})) = \hat{\kappa}(\hat{\vartheta}) - a_*(\vartheta_0) \frac{1}{n} \sum_{j=1}^{n} X_j + o_p(n^{-1/2}),$$

and Theorem 4.1, Remark 4.2 and Lemma 4.1 imply

$$\hat{\kappa}(\hat{\vartheta}) = \tilde{\kappa} + E_P[h'(Y_1(\vartheta_0)) \dot{Y}_1(\vartheta_0)^{\top}](\hat{\vartheta} - \vartheta_0) + o_p(n^{-1/2}),$$

with $\dot{Y}_1(\vartheta_0) = \sum_{r=1}^{\infty} \dot{\delta}_r(\vartheta_0) X_{1-r}$. By Theorem 2.1 we have

$$\tilde{\kappa} = \kappa(\vartheta_0, P) + \frac{1}{n}\sum_{j=1}^{n} h_*(X_j).$$

We arrive at the following result.

THEOREM 5.1. *Suppose assumptions* (5.3) *and* (5.4) *hold and h satisfies Assumption* H [*with* $Y_1(\vartheta_0)$ *playing the role of* $S(\vartheta_0)$]. *Choose* $m = m(n)$ *such that* $m^4/n \to 0$ *and* $\log(n)/m \to 0$. *If* $\hat{\vartheta}$ *is* $n^{1/2}$-*consistent for* $\vartheta_0$, *then*

$$\hat{\tilde{\kappa}}(\hat{\vartheta}, \hat{a}_*(\hat{\vartheta})) = \kappa(\vartheta_0, P) + E_P[h'(Y_1(\vartheta_0))\dot{Y}_1(\vartheta_0)^\top](\hat{\vartheta} - \vartheta_0)$$
$$+ \frac{1}{n}\sum_{j=1}^{n}[h_*(X_j) - a_*(\vartheta_0)X_j] + o_p(n^{-1/2}).$$

Computations are faster if $m$ is small. We may choose $m$ proportional to $(\log n)^{1+\varepsilon}$ with $\varepsilon > 0$.

Let us now show that $\hat{\tilde{\kappa}}(\hat{\vartheta}, \hat{a}_*(\hat{\vartheta}))$ is efficient for $E_P[h(Y_1(\vartheta_0))]$ if $\hat{\vartheta}$ is efficient for $\vartheta_0$. Schick and Wefelmeyer (2002a) give conditions for local asymptotic normality and characterize efficient estimators for differentiable functionals in causal and invertible linear processes. We need only check that the functional $\kappa(\vartheta_0, P) = E_P[h(Y_1(\vartheta_0))]$ is differentiable in an appropriate sense, with efficient influence function equal to the influence function of $\hat{\tilde{\kappa}}(\hat{\vartheta}, \hat{a}_*(\hat{\vartheta}))$.

We assume from now on that $P$ has finite Fisher information $I(P)$ for location; that is, $P$ has an absolutely continuous density $f$ and $I(P) = \int \ell^2 \, dP < \infty$, where $\ell = f'/f$. We also assume that the matrix $V(\vartheta_0) = E_P[\xi_1 \xi_1^\top]$ is positive definite.

Local asymptotic normality and differentiability require a local model. It is introduced in Schick and Wefelmeyer (2002a) as follows. Set

$$G = \left\{g \in L_*(P): \int g \, dP = \int xg(x) \, P(dx) = 0\right\}.$$

For $g$ in $G$ define $P_{n,g}$ by its $P$-density $1 + n^{-1/2} g_n$ with

$$g_n = \overline{g}_n - \int \overline{g}_n \gamma_n^\top \, dP \left(\int \gamma_n \gamma^\top \, dP\right)^{-1} \gamma_n,$$

where $\gamma(x) = (1, x)^\top$ and $\gamma_n(x) = (1, -n^{1/8} \vee x \wedge n^{1/8})^\top$, and

$$\overline{g}_n = \int g \mathbf{1}[|g| \leq n^{1/8}](x - n^{-1/8}y)\varphi(y) \, dy,$$

with $\varphi$ the standard normal density. Set $\vartheta_{n,t} = \vartheta_0 + n^{-1/2}t$ for $t \in \mathbb{R}^d$. The arguments of Theorems 2.2 and 4.1 yield the following result.

THEOREM 5.2. *Suppose assumptions (5.3) and (5.4) hold and h satisfies Assumption H [with $Y_1(\vartheta_0)$ playing the role of $S(\vartheta_0)$]. Then, for each $(t, g) \in \mathbb{R}^d \times G$,*

$$n^{1/2}(\kappa(\vartheta_{n,t}, P_{n,g}) - \kappa(\vartheta_0, P)) \to E_P[h'(Y_1(\vartheta_0))\dot{Y}_1(\vartheta_0)^\top]t + \int h_* g \, dP.$$

Schick and Wefelmeyer [(2002a), Section 5], construct a least dispersed regular estimator $\hat{\vartheta}_*$ for $\vartheta_0$. It is asymptotically linear,

$$\hat{\vartheta}_* = \vartheta_0 + \frac{1}{n} \sum_{j=1}^{n} (V(\vartheta_0)I(P))^{-1} \xi_j \ell(X_j) + o_p(n^{-1/2}).$$

By Theorem 5.1, the substitution estimator $\hat{\hat{\kappa}}(\hat{\vartheta}_*, \hat{a}_*(\hat{\vartheta}_*))$ is also asymptotically linear,

$$\hat{\hat{\kappa}}(\hat{\vartheta}_*, \hat{a}_*(\hat{\vartheta}_*))$$
$$= \kappa(\vartheta_0, P) + \frac{1}{n} \sum_{j=1}^{n} \Big[ E_P[h'(Y_1(\vartheta_0))\dot{Y}_1(\vartheta_0)^\top]$$
$$\times (V(\vartheta_0)I(P))^{-1} \xi_j \ell(X_j) + h_*(X_j) - a_*(\vartheta_0)X_j \Big]$$
$$+ o_p(n^{-1/2}).$$

By the characterization in Schick and Wefelmeyer [(2002a), Section 2], Theorem 5.2 shows that the efficient influence function of $\kappa(\vartheta, P)$ equals the influence function of the substitution estimator $\hat{\hat{\kappa}}(\hat{\vartheta}_*, \hat{a}_*(\hat{\vartheta}_*))$, so that the latter is least dispersed and regular for $\kappa(\vartheta_0, P) = E_P[h(Y_1(\vartheta_0))]$.

**6. Variance reduction in a special case.** We illustrate our results with the autoregressive example considered in the Introduction. Let $Y_0, \ldots, Y_n$ be observations from the AR(1) model $Y_t = \vartheta_0 Y_{t-1} + X_t$ with $|\vartheta_0| < 1$ and independent and identically distributed innovations $X_t$ with distribution $P$, density $f$, mean 0 and finite fourth moment $\mu_4$, where $\mu_k = \int x^k P(dx)$, $k = 2, 3, 4$. We also assume that $P$ has finite Fisher information $I(P) = \int \ell^2 \, dP$ for location, where $\ell = f'/f$.

We want to estimate the stationary variance

$$\sigma^2 = \kappa(\vartheta_0, P) = E[Y_1^2] = E\left[\left(\sum_{s=0}^{\infty} \vartheta_0^s X_s\right)^2\right].$$

Here $h(x) = x^2$. The stationary variance reduces to

$$\sigma^2 = \frac{\mu_2}{1 - \vartheta_0^2}.$$

We consider the following estimators. The empirical estimator of $\sigma^2$ is

$$\hat{\sigma}^2 = \frac{1}{n} \sum_{j=1}^{n} Y_j^2$$

and has influence function

$$\frac{1}{1-\vartheta_0^2}(y^2 - \vartheta_0^2 x^2 - \mu_2).$$

The *improved* empirical estimator of $\sigma^2$ is

$$\hat{\sigma}_*^2 = \frac{1}{n} \sum_{j=1}^{n} \left( Y_j^2 - \frac{\hat{\mu}_3}{(1+\hat{\vartheta}_*)\hat{\mu}_2} Y_j \right),$$

with $\hat{\mu}_k$ as defined in (1.1) and $\hat{\vartheta}_*$ the least squares estimator:

$$\hat{\vartheta}_* = \frac{\sum_{j=1}^{n} Y_{j-1} Y_j}{\sum_{j=1}^{n} Y_{j-1}^2}.$$

The improved empirical estimator has influence function

$$\frac{1}{1-\vartheta_0^2}\left(y^2 - \vartheta_0^2 x^2 - \mu_2 - \frac{\mu_3}{\mu_2}(y - \vartheta_0 x)\right).$$

For these results we refer to Example 2 in Müller, Schick and Wefelmeyer (2001b).

Finally, we write $\hat{\sigma}_\#^2(\hat{\vartheta}) = \hat{\hat{\kappa}}(\hat{\vartheta}, \hat{a}_*(\hat{\vartheta}))$ for our estimator of $\sigma^2$. Suppose that $\hat{\vartheta}$ is asymptotically linear with influence function $w$. Then by Theorem 5.1 our estimator is asymptotically linear with influence function

$$\frac{1}{1-\vartheta_0^2}\left(\frac{2\vartheta_0\mu_2}{1-\vartheta_0^2}w(x,y) + (y - \vartheta_0 x)^2 - \mu_2 - \frac{\mu_3}{\mu_2}(y - \vartheta_0 x)\right).$$

The least squares estimator $\hat{\vartheta}_*$ has influence function

$$w(x,y) = \frac{1-\vartheta_0^2}{\mu_2} x(y - \vartheta_0 x).$$

An efficient estimator $\hat{\vartheta}_\#$ has influence function

$$w(x,y) = -\frac{1-\vartheta_0^2}{\mu_2 I(P)} x\ell(y - \vartheta_0 x).$$

If we use an efficient estimator $\hat{\vartheta}_\#$, then the estimator $\hat{\sigma}_\#^2(\hat{\theta}_\#)$ is efficient by Section 5. In the particular case of estimating moments, simpler efficient estimators are given in Section 6 of Schick and Wefelmeyer (2002a). In particular, a simpler efficient estimator of $\sigma^2$ is

$$\frac{\hat{\mu}_2^*}{1-\hat{\vartheta}_\#^2} \qquad \text{with } \hat{\mu}_2^* = \hat{\mu}_2 - \frac{\hat{\mu}_3}{\hat{\mu}_2}\hat{\mu}_1.$$

The estimator is obtained by replacing $\mu_2$ and $\vartheta_0$ in $\sigma^2 = \mu_2/(1-\vartheta_0^2)$ by efficient estimators. The efficient estimator $\hat{\mu}_2^*$ of $\mu_2$ uses the constraint $\mu_1 = 0$. Of course, both efficient estimators for $\sigma^2$ are stochastically equivalent. This can be seen directly by simplifying $\hat{\sigma}_\#^2(\hat{\vartheta}_\#)$. More generally, $\hat{\mu}_2^*/(1-\hat{\vartheta}^2)$ is stochastically equivalent to $\hat{\sigma}_\#^2(\hat{\vartheta})$ for any $n^{1/2}$-consistent estimator $\hat{\vartheta}$ of $\vartheta_0$.

Next we determine the asymptotic variances of these estimators. The empirical estimator $\hat{\sigma}^2$ has asymptotic variance
$$\frac{1}{(1-\vartheta_0^2)^2}\left(\mu_4 - \mu_2^2 + 4\mu_2^2\frac{\vartheta_0^2}{1-\vartheta_0^2}\right).$$

The *improved* empirical estimator $\hat{\sigma}_*^2$ has asymptotic variance
$$\frac{1}{(1-\vartheta_0^2)^2}\left(\mu_4 - \mu_2^2 + 4\mu_2^2\frac{\vartheta_0^2}{1-\vartheta_0^2} - \frac{\mu_3^2}{\mu_2}\right).$$

One calculates that the estimators $\hat{\sigma}_\#^2(\hat{\theta}_*)$ and $\hat{\mu}_2^*/(1-\hat{\vartheta}_*^2)$ have the same asymptotic variance. Finally, the efficient estimators $\hat{\sigma}_\#^2(\hat{\theta}_\#)$ and $\hat{\mu}_2^*/(1-\hat{\vartheta}_\#^2)$ have asymptotic variance
$$\frac{1}{(1-\vartheta_0^2)^2}\left(\mu_4 - \mu_2^2 + \frac{4\mu_2\vartheta_0^2}{I(P)(1-\vartheta_0^2)} - \frac{\mu_3^2}{\mu_2}\right).$$

The relative asymptotic variance increase of the empirical estimator $\hat{\sigma}^2$ over the efficient estimator is
$$\frac{I(P)(1-\vartheta_0^2)\mu_3^2/\mu_2 + 4\vartheta_0^2\mu_2(\mu_2 I(P) - 1)}{I(P)(1-\vartheta_0^2)(\mu_4 - \mu_2^2 - \mu_3^2/\mu_2) + 4\vartheta_0^2\mu_2}.$$

For the improved empirical estimator $\hat{\sigma}_*^2$, the relative asymptotic variance increase is
$$\frac{4\vartheta_0^2\mu_2(\mu_2 I(P) - 1)}{I(P)(1-\vartheta_0^2)(\mu_4 - \mu_2^2 - \mu_3^2/\mu_2) + 4\vartheta_0^2\mu_2}.$$

These estimators are efficient for values of $\vartheta_0$ and $P$ for which the corresponding ratios are 0. The second ratio is 0 if and only if $\vartheta_0 = 0$ or $\mu_2 I(P) = 1$. The latter happens if and only if $P$ is normal. Thus the improved empirical estimator $\hat{\sigma}_*^2$ is efficient if and only if $\vartheta_0 = 0$ or $P$ is normal. The first ratio is 0 if and only if $\mu_3 = 0$ and also $\vartheta_0 = 0$ and $\mu_2 I(P) = 1$. Thus the empirical estimator $\hat{\sigma}^2$ is efficient if $P$ is the normal distribution. For other distributions, it is efficient if and only if $\vartheta_0 = 0$ and $\mu_3 = 0$. The two ratios are the same if and only if $\mu_3 = 0$, which is the case for symmetric $P$.

If $\vartheta_0$ is close to 1, both ratios are close to $\mu_2 I(P) - 1$. Note that $\mu_2 I(P) - 1$ is the relative variance increase of the sample mean versus the efficient estimator in the location model generated by $P$. It is well known that $\mu_2 I(P) - 1$ can be large if $P$ is not normal.

**Acknowledgments.** We thank the referees for helpful comments and suggestions.

Department of Mathematical Sciences  
Binghamton University  
Binghamton, New York 13902-6000  
USA  
e-mail: anton@math.binghamton.edu

Mathematisches Institut  
Universität zu Köln  
Weyertal 86-90  
50931 Köln  
Germany  
e-mail: wefelm@math.uni-koeln.de